\documentclass[11pt]{amsart}
\hoffset= - 0.75 in

\usepackage{amsfonts, amsmath, amssymb, amsthm}
\usepackage{latexsym, graphicx, graphics, color}
\newtheorem{thm}{Theorem}[section]
\newtheorem{move}{Move}

\newtheorem{question}[thm]{Question}

\newtheorem{defn}[thm]{Definition}

\newtheorem{lem}[thm]{Lemma}
\newtheorem{cor}[thm]{Corollary}
\newtheorem{prop}[thm]{Proposition}

\newcommand{\ZZ}{\mathbb{Z}}
\newcommand{\CC}{\mathbb{C}}

\begin{document}

\title{Rook poset equivalence of Ferrers boards}
\author{Mike Develin}
\address{Mike Develin, American Institute of Mathematics, 360 Portage 
Ave., Palo Alto, CA 94306-2244, USA}
\date{\today}
\email{develin@post.harvard.edu}

\begin{abstract}

A natural construction due to K. Ding yields Schubert varieties from 
Ferrers boards. The poset structure of the Schubert cells in these 
varieties is equal to the poset of maximal rook placements on the Ferrers 
board under the Bruhat order. We determine when two Ferrers boards have 
isomorphic rook posets. Equivalently, we give an exact categorization of 
when two Ding Schubert varieties have identical Schubert cell structures. 
This also produces a complete classification of isomorphism types of lower 
intervals of 312-avoiding permutations in the Bruhat order.

\end{abstract}

\maketitle

\section{Introduction and Definitions}

We begin by introducing our main object of study.

\begin{figure} 
\begin{center}
\includegraphics{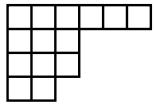}
\end{center}
\caption{\label{ferrers}
The Ferrers board corresponding to the partition (2, 3, 3, 6).
}
\end{figure}

\begin{defn}
Let $\lambda$ be a partition $0\le \lambda_1\le\cdots\le\lambda_n$, where 
the $\lambda_i$'s are positive integers (note 
that our partitions will always be written with the parts in weakly 
increasing order.) The \textbf{Ferrers board} $B_\lambda$ is the grid 
consisting of $\lambda_i$ left-justified squares in the $i$-th row from 
the bottom, as in 
Figure~\ref{ferrers}.
\end{defn}

A classical combinatorial problem is that of rook-equivalence~\cite{Kapl}: 
when is it the case that, for each $k$, the number of ways to place $k$ 
non-attacking rooks on each of two Ferrers boards is the same? This 
problem was originally solved by Foata and Sch\"utzenberger~\cite{FS}, and 
later reproved by Goldman, Joichi, and White~\cite{GJW}, who showed that 
$B_\lambda$ and $B_\mu$ are rook-equivalent if and only if the multisets 
$\{\lambda_i-i\}_{i=1}^n$ and $\{\mu_i-i\}_{i=1}^n$ are the same. We give 
these sequences of numbers a name.

\begin{defn}
Let $\lambda$ be a partition. The \textbf{GJW sequence} of $\lambda$ is the 
sequence $(\lambda_i-i)_{i=1}^n$. Note that the GJW sequence cannot decrease by more than 1.
\end{defn}

Throughout this paper, we will consider only partitions $\lambda$ which 
contain a staircase, i.e. such that for all $1\le i\le n$, $\lambda_i \ge 
i$. This is equivalent to the GJW-partition containing only nonnegative 
integers; these partitions are the ones which can actually accommodate $n$ 
rooks. To such a board, we can associate a Schubert variety $X_\lambda$ 
inside a certain partial flag manifold,
consisting of all flags $0\subset V_1\subset V_2\subset \cdots 
\subset V_n\subset \CC^N$ such that dim($V_i$) = i and $V_i$ is contained 
inside the subspace of $\CC^n$ generated by the first $\lambda_i$ basis 
vectors. This construction is due to Ding~\cite{Ding}, who applied work of 
Garsia and Remmel~\cite{GR} to compute the Poincare series of these 
varieties. In earlier work with J. Martin and V. Reiner, we 
classified these varieties up to isomorphism.

\begin{thm}\cite{DMR}\label{DMR-thm}
Two Schubert varieties $X_\lambda$ and $X_\mu$ are isomorphic if and only if upon 
removing 0's from their GJW sequences, the remaining multisets of contiguous blocks 
are identical. Indeed, this is the only way they can have isomorphic cohomology 
rings over $\ZZ$.
\end{thm}

0's in the GJW sequence correspond to subspaces in the flag where $V_i$ 
is 
completely specified; the decomposition into blocks given in Theorem~\ref{DMR-thm} 
therefore expresses the original Schubert variety as a product of 
constituents given by the blocks. The content of the theorem is that this is the 
{\bf only} way that two boards can yield isomorphic Schubert varieties, or even 
ones with isomorphic cohomology rings.

We can say a lot about these Ding Schubert varieties simply from the combinatorics 
of the Ferrers board. The following poset is the key ingredient.

\begin{defn}
To each maximal rook placement $x$ on a Ferrers board $\lambda$, we associate a 
sequence 
$[x_i]_{i=1}^n$ where $x_i$ is the position (i.e. column index) of the 
rook in the $i$-th row. The 
\textbf{rook poset} $P_\lambda$ has elements given by the maximal rook placements, 
with $x\le y$ if for each $j$, the set $\{x_1, \ldots, x_j\}$, placed in increasing 
order, is positionwise less than or equal to the set $\{y_1, \ldots, y_j\}$ placed 
in 
increasing order.
\end{defn}

For example, consider the Ferrers board given by the partition $(3, 3, 5, 
6, 6)$. The rook placement $[2, 1, 5, 3, 4]$ lies below the rook placement 
$[3, 1, 5, 2, 6]$ in the rook poset: sorting each initial segment, we find 
that $2 \le 3$, $\{1, 2\}\le \{1, 3\}$, $\{1, 2, 5\}\le \{1, 3, 5\}$, 
$\{1, 2, 3, 5\}\le \{1, 2, 3, 5\}$, and $\{1, 2, 3, 4, 5\}\le \{1, 2, 3, 
5, 6\}$. Some statements about the rook poset are immediate:

\begin{prop} 
The rook poset has a maximal element $\hat{1}$, given by placing each rook as far 
to the 
right as it will go, starting from the bottom, and a minimal element $\hat{0}$, 
given by 
placing the rook in row $i$ in position $i$.
\end{prop}

For instance, the maximal element of (3, 3, 5, 6, 6) is [3, 2, 5, 6, 4]. 
Each initial segment is easily seen to be maximal among all initial 
segments of rook placements, and hence this element is maximal. The 
minimal element of (3, 3, 5, 6, 6) is [1, 2, 3, 4, 5], again trivial to 
verify.

Our goal is to investigate when two such rook posets are isomorphic. The following 
definition and proposition provide one such case.

\begin{defn}

Let $\lambda$ be a partition with $\lambda_n = n+1$. Define the \textbf{conjugate} 
of $\lambda$ as follows: take the Ferrers board of $\lambda$, add a row of $n+1$ 
boxes at the top, reflect the diagram (which now has width and height both equal to 
$n+1$) across the upper-left to lower-right diagonal, and remove the top row of 
$n+1$ boxes. The resulting Ferrers board is the conjugate of $\lambda$. We 
similarly call two GJW sequences ending in 1 conjugate if their underlying 
partitions are. 

\end{defn}

\begin{figure} 
\begin{center}
\includegraphics{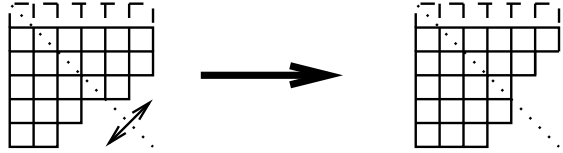}
\end{center}
\caption{\label{conjugatefig}
Conjugate partitions. The one on the left has $\lambda=(2, 3, 5, 6, 6)$, 
and it is conjugate to the one on the right, which has $\lambda = (3, 4, 
4, 5, 6)$.}
\end{figure}

The following proposition is relatively straightforward.

\begin{prop}\label{conjugate}
The rook posets of $\lambda$ and of its conjugate are isomorphic. The 
isomorphism is given 
by taking a maximal rook placement on $\lambda$, extending it to a maximal rook placement on $(\lambda, n+1)$ in the only way possible, conjugating the rook placement, and removing the top row.
\end{prop}

Less obvious is the natural grading placed on the rook poset.

\begin{prop}\cite{GR}
Given a maximal rook placement $x$, complete the its sequence to a permutation 
$\sigma_x$ of 
$\{1, 2, \ldots, \lambda_n\}$ by sending each of $n+1, \ldots, 
\lambda_n$ in turn to the smallest possible value. Define $r(x)$ to be the number 
of inversions in this permutation, i.e. the number of $i<j$ with $\sigma_x(i) > 
\sigma_x(j)$. Then $r(x)$ defines a grading on the rook poset.
\end{prop}

\begin{figure}
\begin{center}
\includegraphics{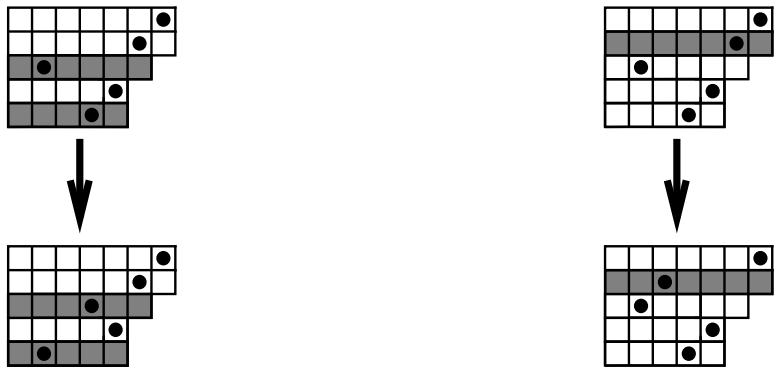}
\end{center}
\caption{\label{moves}
A switch move (left) and push move (right).
}
\end{figure}

Ding~\cite{Ding2} showed that the rank sizes under this grading yield the Poincar\'e 
series of the corresponding Schubert variety. In addition, he showed that an 
individual rook placement corresponds to a Schubert cell, whose dimension is equal 
to the rank of the rook placement; the closure of this Schubert cell 
consists of the union of all of the Schubert cells corresponding to rook 
placements lying below it in the rook poset.
The covering relations are given by single inversions, which amount 
to performing one of the following two operations to obtain a maximal rook 
placement immediately below a given maximal rook placement (see 
Figure~\ref{moves}.)

\begin{move}\label{switch}
If rooks lie in positions $(i, j)$ and $(k, l)$ with $i<k$ and $j>l$, and 
the 
rectangle defined by these two squares empty aside from the two rooks at its 
corners, switch these rooks out in favor of $(i, l)$ and $(j, k)$. We call this a 
\textbf{switch move} on rook $i$. \end{move}
\begin{move}\label{push}
If a rook lies in position $(i, j)$, column $k$ is empty with $k<j$, and for all 
$k<r<j$, there is a rook in column $r$ below row $i$, then replace the rook at $(i, 
j)$ with one at $(i, k)$. We call this a \textbf{push move} on rook $i$.
\end{move}

If the last element of the GJW sequence is a 1 or a 0, the rook poset is 
isomorphic to the interval $[\sigma_{\hat{0}}, \sigma_{\hat{1}}]$ in the 
Bruhat order on the symmetric group $S_{\lambda_n}$; if it is not, then it 
is isomorphic to a lower interval in the Bruhat order of the associated 
partial flag manifold. The permutation $\sigma_{\hat{0}}$ is of course the 
identity, so 
investigating the structure of rook posets is tantamount to investigating 
lower intervals in the Bruhat order. Focusing on the first case above, the 
permutations $\sigma_{\hat{1}}$ are precisely those which are 
312-avoiding, an important class of permutations. A description of when 
two Ferrers boards have isomorphic rook posets therefore produces as a 
corollary a description of when the lower intervals of 312-avoiding 
permutations in the Bruhat order are isomorphic.

In this paper, we establish necessary and sufficient conditions for two 
rook posets to be isomorphic. This establishes a notion of equivalence of 
Ferrers boards stronger than ordinary rook-equivalence and weaker than the 
finer rook-equivalence given by isomorphism of the corresponding Schubert 
varieties. Aside from being of intrinsic combinatorial interest, this has 
two other main consequences: it gives us geometric information about the 
Schubert cell structures of the Ding varieties, and it gives us an exact 
classification of lower intervals of 312-avoiding permutations in the 
Bruhat order up to isomorphism.

\section{The Main Theorem}

After extensive experimentation, Martin, Reiner, and Wagner~\cite{MRW} 
formulated the conjecture which is the main theorem of this paper.

\begin{thm}\label{mainthm}
Let $\lambda$ and $\mu$ be two partitions. Then the corresponding Ferrers boards 
have isomorphic rook posets if and only if when we break the GJW sequences of 
$\lambda$ and $\mu$ after every 1, and throw out all zeroes, the resulting 
multisets of contiguous blocks are identical up to conjugation.
\end{thm}

For instance, the GJW sequences 22321010032 and 13222132 correspond to Ferrers 
boards with isomorphic posets: the first sequence breaks into 22321, 1, and 32 while the second breaks into 1, 32221, and 32. 22321 and 32221 turn out to be conjugate, while the other two blocks are identical. Note that only one block can not end in a 1, the last one; since all 0's must follow 1's (because the sequence cannot decrease by more than one), all other blocks end in 1's.

It is relatively straightforward to show one direction, namely that if two 
GJW sequences produce identical block multisets up to conjugation, then 
their rook posets are equivalent. The remainder of this section will prove 
this by checking each operation involved.

\begin{lem}\label{productlem}
Let $\lambda = (\lambda^1, 0, \lambda^2)$ be a partition. Then the rook 
poset of $\lambda$ is isomorphic to the product of the rook posets of 
$\lambda^1$ and $\lambda^2$.
\end{lem}
\begin{proof}
Because of the existence of the 0 row, a maximal rook placement on 
$\lambda$ simply consists of the union of a maximal rook placement on 
$(\lambda^1, 0)$ and a maximal rook placement on $\lambda^2$ whose rooks 
all lie to the right of the rightmost square in $(\lambda^1, 0)$. The two 
parts of the maximal rook placement do not interact in any way, and the 
statement follows easily; the rook posets of $(\lambda^1, 0)$ and 
$\lambda^1$ are easily seen to be isomorphic.
\end{proof}

Indeed, not only the rook posets, but also the Schubert varieties and 
cohomology rings are isomorphic. Less trivial is the fact that we can 
break the blocks after 1's.

\begin{lem}\label{ins0lem}
Let $\lambda = (\lambda^1, 1, \lambda^2)$ be a partition. Then $\lambda$ and $\lambda' = (\lambda^1, 1, 0, \lambda^2)$ have isomorphic rook posets.
\end{lem}

\begin{figure}
\begin{center}
\includegraphics{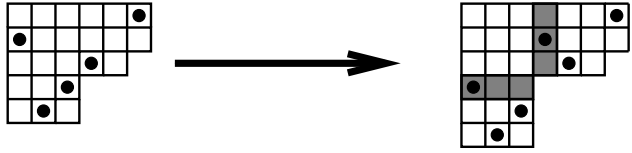}
\end{center}
\caption{\label{insertion}
Isomorphism of rook posets under the insertion of a 0 after a 1 in the GJW 
sequence. A row is inserted to create the 0, and its rook is uniquely 
specified; the rook formerly in that column is moved to an inserted column immediately to the 
right of the inserted row. The GJW sequence of the new diagram (here 
$(2, 1, 0, 2, 2, 1)$) is the same 
as that of the old one (here $(2, 1, 2, 2, 1)$) with an extra 0 inserted 
after a 1.}
\end{figure}

\begin{proof}
Let the number of parts in $\lambda^1$ be $n-2$, so that the inserted row 
is in position $n$ and has length $n$. We describe the isomorphism. Given 
a maximal rook placement $x = (x_1, \ldots, x_m)$ on $\lambda$, we define 
a maximal rook placement $x'$ on $\lambda'$ as follows; see 
Figure~\ref{insertion}.
\begin{itemize}
\item For $1\le i\le n-1$, $x_i' = x_i$.
\item Let $x_n'$ be the unique integer from $[n]$ which is not one of $(x_1, \ldots, x_{n-1})$. 
\item For $i>n$, $x_i' = x_{i-1} + 1$, unless $x_{i-1} = x_n'$ (i.e. is the missing column from $[n]$), in which case $x_i' = n+1$. 
\end{itemize}
It is straightforward to verify that this procedure yields a maximal rook placement. The fact that this operation preserves the rook poset structure is similarly easy to check, containing a couple of cases, and is left to the reader. Finally, the map is clearly a bijection, since we can recover $x$ from $x'$ by removing the $n$-th row, shortening the $(n+1)$-st and higher rows by one box (on the left; these boxes must be empty), and moving any rook formerly in the $(n+1)$-st column to the column formerly occupied by the excised $n$-th row rook.
\end{proof}

By inserting 0's after all the 1's via Lemma~\ref{ins0lem}, and applying 
Lemma~\ref{productlem}, we conclude that the rook poset of a partition is 
equal to the product of the rook posets of its blocks (as defined in 
Theorem~\ref{mainthm}.) By Proposition~\ref{conjugate}, conjugate 
partitions have isomorphic rook posets, which completes the proof of the 
easy direction of the theorem.

\section{The Hard Direction}

The hard direction of Theorem~\ref{mainthm} is proving that 
two partitions which do not decompose identically into 
blocks, up to conjugation, actually have non-isomorphic 
posets. Showing that two posets are not isomorphic is generally rather ad-hoc; the 
Garsia-Remmel and Ding results determine when the gradings of these posets 
yield different rank sizes (one easy way to show that two posets are not 
isomorphic), namely when their GJW multisets are different, but the 
conditions of Theorem~\ref{mainthm} are much more restrictive. 

Our approach is to use a case-by-case analysis to rule out all other 
isomorphisms, reconstructing the GJW sequence from the rook poset up to 
the equivalence of Theorem~\ref{mainthm}. This approach will also yield 
intrinsic information about the poset (and therefore the Ding Schubert 
variety) which can be read off from the GJW sequence. At times, we will 
use an induction on the number of boxes in the diagram, assuming that we 
can distinguish the rook posets of smaller diagrams which do not satisfy 
the hypotheses of Theorem~\ref{mainthm}.

Because of Lemma~\ref{ins0lem}, we will assume that there is exactly one 0 
after every 1. This means that the x-coordinates of the rooks in distinct 
blocks lie in distinct intervals. Our main objects in the investigation 
will be the \textbf{coatoms}, the elements which have rank one less than 
the maximal element $\hat{1}$.

\begin{prop}\label{coatomrow}
There is a natural bijection between the non-zero elements of the GJW sequence and the coatoms of the corresponding rook poset.
\end{prop}
\begin{proof}

The maximal element $\hat{1}$ has its rooks as far to the right as possible. For 
each row, we need to enumerate the switch moves on the rook in that row, as well as 
the push moves on it. We claim that for each row $i$ corresponding to a non-zero 
GJW element there is exactly one move, and for each row corresponding to a 0 in the 
GJW sequence, there are no moves.

Suppose that row $i$ corresponds to a 0 in the GJW sequence. Then $\lambda_i = i$, 
and each of the rooks in the first $i$ rows lies in the first $i$ columns. 
Therefore, each of the first $i$ columns contains a rook in or below row $i$, and 
so none can contain a rook above row $i$. Since the rook in row $i$ is in one of 
these columns, there are no rooks above and to the left of it, so there are no 
switch moves on rook $i$. Similarly, every column to the left of the 
rook in row $i$'s column is occupied by a rook below it, so there are no push moves 
on rook $i$.

Next, suppose that row $i$ does not correspond to a 0 in the GJW sequence. Suppose 
that there are two rooks above and to the left of it, corresponding to two switch 
moves. Since both rectangles corresponding to the switch moves must be empty, this 
means that we have rooks at $(i, j), (r, s)$, and $(x, y)$ with $i < r < x$ and $s 
< y < j$. This is impossible, since in $\hat{1}$ the rooks are placed as far to the 
right in each row as possible, bottom up; since $r < x$, and there are no rooks in 
column $y$ below row $x$ and hence none below row $r$, we could have placed a rook 
at $(r, y)$. Similarly, if there is a switch move switching $(i, j)$ with $(x, y)$, 
then there cannot be a push move sending $(i, j)$ to $(i, k)$: for $k<y$, there is 
a rook in an intervening column above row $i$, while for $k>y$, this means that 
column $k$ is empty in $\hat{1}$, so we would have placed $(x, y)$ at $(x, k)$. So 
there is at most one move associated to the rook in row $i$.

Suppose there is no switch move on the rook in row $i$; let this rook be in column 
$j$. Because row $i$ does not correspond to a 0 in the GJW sequence, we have 
$\lambda_i \ge i+1$. Since there is no switch move on the rook in row 
$i$, there are no rooks above row $i$ in columns $1, \ldots, j-1$. Since we must have placed the rook as far to the right as possible, there must be rooks below row $i$ in columns $j+1, \ldots, i+1$, a total of $j-i+1$ rooks. This means that, since there are $i-1$ rooks below row $i$, at most $(i-1) - (j-i+1) = j-2$ of them lie in columns $1, \ldots, j-1$. Therefore, there is an empty column among the first $j-1$, so there is a push move on $(i, j)$, completing the proof of the lemma.
\end{proof}

Proposition~\ref{coatomrow} establishes a bijection between the coatoms and the nonzero elements of the GJW sequence (or the rows); this yields a natural labeling of the coatoms with their rows. We can use the GJW sequence to easily determine whether a coatom is a switch move or a push move.

\begin{prop}\label{whenswitch}
A coatom $i$ is a push move if and only if there is no later element of the GJW sequence which is lower than element $i$. If there is such an element, the row of the first such element is the other row of the switch move.
\end{prop}
\begin{proof}

Being a switch move is equivalent to there being a rook above and to the left of 
the rook in the $i$-th row. This happens if and only if, when we reach some row 
$j>i$, the columns $\lambda_i+1, \ldots, \lambda_j$ are all already filled with 
rooks. Since we are looking for the first occurrence of this phenomenon, we can 
assume that all rooks in rows $i+1, \ldots, j-1$ are to the right of the rook in 
row $i$, so this happens exactly when this number, $j-i-1$, is equal for the 
first time to the number 
of columns to fill, $\lambda_j - \lambda_i$. This is equivalent to the statement 
about GJW sequence elements.

\end{proof}

In practice, we will use this proposition in reverse: once we have determined which 
row is the other row of a given switch move, this tells us something about the GJW 
sequence. The next step is to relate the GJW sequence to differences 
between the structure lying underneath the coatoms.

\begin{defn}

Let $I$ be a set of coatoms. Define the \textbf{poset generated by I} to 
be 
$$X_I = \{x\not\le j \text{ for all }j\notin I\}.$$ Define $X_{=I}$ to be 
the subset of $X_I$ for which $x\le i$ for all $i\in I$.

\end{defn}

In other words, $X_I$ is the portion of the poset which does not lie below any 
coatoms not in $I$, while $X_{=I}$ is the portion of $X_I$ which lies below every 
coatom in $I$. Our goal is to use the structure of various $X_I$'s to recover the 
GJW sequence 
(up to the equivalences given by Theorem~\ref{mainthm}.) A key simplifying lemma 
allows us to easily check membership in $X_I$.

\begin{lem}\label{belowcoat}
Suppose the move given by coatom $i$ sends the rook at $(i, j)$ to $(i, k)$. Then to check if a maximal rook placement $x$ lies below $i$, we need only check whether 
$$\{x_1, \ldots, x_i\}\le \{i_1, \ldots, i_{i-1}, i_i = k\}.$$
\end{lem}

\begin{proof}
For $r<i$, $\{i_1, \ldots, i_r\} = \{\hat{1}_1, \ldots, \hat{1}_r\}$, so $\{x_1, \ldots, x_r\}\le \{i_1, \ldots, i_r\}$ since $\hat{1}$ is the maximal element of the rook poset. Suppose first that $i$ is a push move. Then for $r>i$, the rooks in rows $i+1, \ldots, r$ are all as far right as possible in the maximal rook placement given by coatom $i$, so if $\{x_1, \ldots, x_i\}\le \{i_1, \ldots, i_i\}$, then since we have $x_s\le i_s = \hat{1}_s$ for $i<s\le r$, we will have $\{x_1, \ldots, x_r\}\le \{i_1, \ldots, i_r\}$.

On the other hand, suppose $i$ is a switch move with row $a$. For $i<r<a$, the 
above logic again suffices. For $r\ge a$, $\{i_1, \ldots, i_r\} = \{\hat{1}_1, 
\ldots, \hat{1}_r\}$, since the elements are all the same except for $i_i$ and 
$i_a$, which are switched. So we must have $\{x_1, \ldots, x_r\} \le \{i_1, \ldots, 
i_r\}$.

\end{proof}

\begin{defn}

Two coatoms $i$ and $j$ are \textbf{entangled} if there exist two distinct elements 
of corank two in $X_{=\{i, j\}}$. The \textbf{entanglement graph} of $\lambda$ has 
vertices equal to the coatoms, with an edge between $i$ and $j$ if the two coatoms 
are entangled.

\end{defn}

In other words, two coatoms are entangled if there exist two elements which are 
covered by both of them but by no other coatoms. The following proposition 
establishes when two coatoms are entangled.

\begin{prop}
Two coatoms are entangled precisely when there is a switch move between the corresponding rows in $\hat{1}$.
\end{prop}

\begin{proof}
Let the coatoms be $i$ and $x$ with $i<x$. First, suppose there is a switch move between the rows. Then there are rooks in $\hat{1}$ at $(i, j)$ and $(x, y)$; in $i$ those rooks are replaced by $(i, y)$ and $(x, j)$, while in $x$ $(i, j)$ remains fixed and $(x, y)$ is replaced by $(x, z)$ for some $z<y$, possibly with $(a, z)$ replaced by $(a, y)$ for $a>x$ (if $x$ is a switch move.) We claim that there are two maximal rook placements in $X_{={i, j}}$:

\begin{enumerate}
\item Identical to $\hat{1}$ except with $(i, y)$ and $(x, z)$ [and $(a, j)$] instead of the rooks in those rows.
\item Identical to $\hat{1}$ except with $(i, z)$ and $(x, j)$ [and $(a, y)$] instead of the rooks in those rows.
\end{enumerate}

It is easy to verify that these lie below only the coatoms $i$ and $j$, 
using Lemma~\ref{belowcoat} and the fact that no other rook lies in the rectangle with corners $(i, j)$ and $(x, y)$ (or the rectangle with corners $(x, y)$ and $(a, z)$ if applicable.)

Conversely, suppose that there is no switch move between the rows. In almost all cases, this means that the coatoms $i$ and $x$ involve completely unrelated rows; the only element covered by both of them is the one where $i$ and $x$ are both done. The only case which is not straightforward is where $i$ and $x$ are both switch moves with a common third row, so that in $\hat{1}$ there are rooks at $(i, j)$, $(x, y)$, and $(a, b)$, with $i < x < a$ and $b < j < y$. In this case, the only element which lies beneath both $i$ and $x$ is the one where these rooks are replaced by $(i, b)$, $(x, j)$, and $(a, y)$ and the rest of the placement remains unchanged. Again, this is easy to check by application of Lemma~\ref{belowcoat}.
\end{proof}

\begin{figure}
\begin{center}
\includegraphics{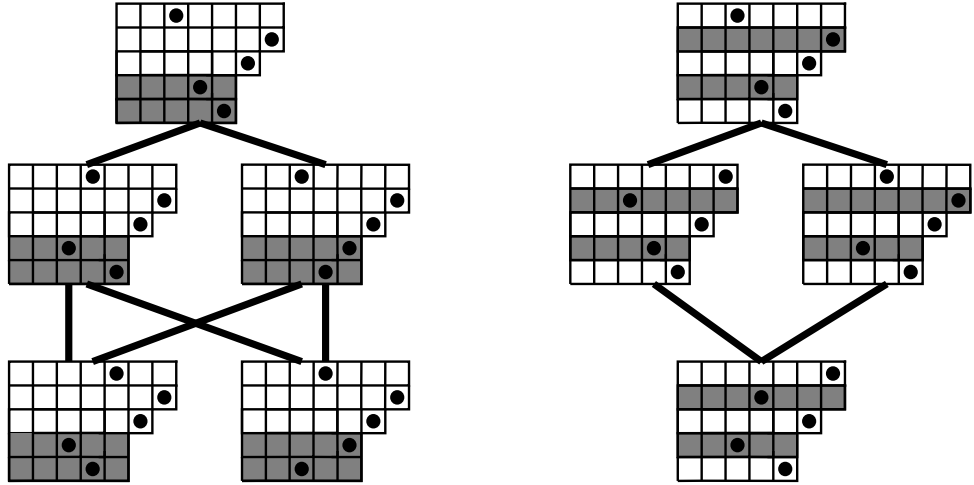}
\end{center}
\caption{\label{entangle}
The relevant part of the rook poset for entangled coatoms (left) and 
unentangled coatoms (right).} 
\end{figure}

\begin{cor}\label{tree}
The entanglement graph of $\lambda$ is a forest.
\end{cor}
\begin{proof}
Every coatom has at most one edge connecting it to a coatom which corresponds to a 
higher row. If the entanglement graph contained a cycle, the lowest coatom in it 
would violate this criterion.
\end{proof}

Each component of the forest has a natural root: the row which is highest among the 
rows comprising that component. Placing these roots at the top, beneath each 
element of the tree, there is a canonical way to order the children: in increasing 
order of their row number. This \textbf{ordered entanglement graph} is easily drawn 
from the diagram.

\begin{prop}
The ordered entanglement graph of $\lambda$ can be drawn from the maximal rook 
placement $\hat{1}$ as follows: draw a segment connecting any two rooks connected 
by a switch move, and turn the diagram forty-five degrees clockwise, ignoring all segments ending 
in a 0 row.
\end{prop}

\begin{proof}
The proof is straightforward. If two segments crossed, this would imply the 
existence of multiple switch moves on some rook, which cannot happen by 
Proposition~\ref{coatomrow}. If one child lies above another, it also must lie to 
the right, and hence must have higher coordinate sum (which is the left-to-right 
ordering given by the forty-five degree turn.)
\end{proof}

\begin{figure}
\begin{center}
\includegraphics{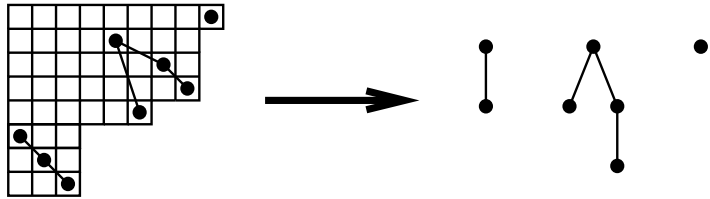}
\end{center}
\caption{\label{entgraph}
How to construct the ordered entanglement graph of a diagram.}
\end{figure}

By definition, we can obtain the entanglement graph from the rook poset. 
We now set upon the task of recovering the ordered entanglement graph from 
the entanglement graph, and the GJW sequence from the ordered entanglement 
graph, up to the equivalence of Theorem~\ref{mainthm}. Each connected 
component of the entanglement graph corresponds to a contiguous set of 
rows; our first goal is to 
recover the blocks, each of which consists of the union of some subset of 
the connected components.

If every block ends in a 1; then the connected components of the entanglement graph 
of $\lambda$ are precisely the blocks (by Proposition~\ref{whenswitch}). The only 
possible problem is introduced by the last block, which may not end in a 1 and thus 
may consist of several connected components of the entanglement graph. Suppose that 
a connected component lies in the last block, which does not end in a 1. Its 
uppermost row is then a number greater than 1, with the property that there is no 
number smaller than it which follows it in the GJW sequence. By 
Proposition~\ref{whenswitch}, it must be a push move. We can determine even more 
about it.

\begin{prop}\label{whenpush}
Let $i$ be a coatom. Then $X_i$ is a chain whose length is 1 if $i$ is a switch move, and the $i$-th element of the GJW sequence if $i$ is a push move.
\end{prop}
\begin{proof}

Suppose first that $i$ is a switch move switching $(i, j)$ and $(x, y)$ to 
$(i, y)$ and $(x, j)$. If row $x$ corresponds to a 0 in the GJW sequence, 
then it is easy to see that no more moves can be performed on the rook in 
row $i$, and performing a move on any other rook creates a position below 
the corresponding coatom. So we reduce to the case where coatom $x$ 
exists. Suppose that coatom $x$ sends $(x, y)$ to $(x, z)$. We need to 
show that any element of the rook poset lying strictly below $i$ lies 
below some other coatom; here again Lemma~\ref{belowcoat} greatly 
simplifies 
the calculations which we omit. If we change the position of any rook 
other than $i$ and $x$, it is easy to see that we will be below the coatom 
corresponding to that row. Suppose we move rook $i$ again, via either a 
switch move or a push move, to column $r<y$; since $(i, r)$ and $(x, j)$ 
are a valid placement with the rest of $\hat{1}$, $(i, j)$ and $(x, r)$ 
must also be, which implies that $z = r$. The new position, with rooks at 
$(i, z)$ and $(x, j)$, will then be below move $x$.

Suppose we change the position of rook $x$ via a switch move or push move. We claim that every move which can be done on $(x, j)$ in the $i$-placement can also be done on $(x, y)$ in the $\hat{1}$-placement. The only other possibilities are a push move to a column between $y$ and $j$ or a switch move with a rook between columns $y$ and $j$, but both are impossible since all columns between $y$ and $j$ must be filled below $x$ so that we place rook $x$ at $(x, y)$. So the only move on row $x$ in the $i$-placement sends $(x, y)$ to $(x, z)$, and by Lemma~\ref{belowcoat} it is immediate that after this move the resulting element of the rook poset lies below coatom $x$.

Next, suppose $i$ is a push move sending $(i, j)$ to $(i, y)$. By 
Lemma~\ref{belowcoat}, it is easy to see that the elements in $X_i$ are 
precisely those which change only the rook in row $i$, as moving any other rook left puts the resulting position below that coatom, and there are no rooks above and to the left of $(i, y)$ to do a switch move, since such a rook would be above and to the left of $(i, j)$, and there are no switch moves on $(i, j)$ in $\hat{1}$. The poset structure of these elements is clearly a chain, and their number is equal to the number of possible rooks $(i, k)$, $k\le j$, which are compatible with the rest of $\hat{1}$; this is just the number of empty columns to the left of $(i, j)$. Since no rooks are above and to the right of $(i, j)$, the number of occupied columns is precisely the number of rooks below and to the left of $(i, j)$. The number of rooks below $(i, j)$ is $i-1$, and the number of rooks below and to the right of $(i, j)$ is $\lambda_i - j$, so the number of rooks below and to the left of $(i, j)$ is $(i-1) - (\lambda_i - j) = i + j - 1 - \lambda_i$. The number of columns to the left of $(i, j)$ is of course $j-1$, so the number of unoccupied columns is $(j-1) - (i + j - 1 - \lambda_i) = \lambda_i - i$, the $i$-th element of the GJW sequence of $\lambda$ as desired.
\end{proof}

So the push moves are precisely the coatoms for which $X_i$ is nontrivial. This 
yields a way to determine the blocks: if a connected component of the entanglement 
graph of $\lambda$ has $|X_i| = 1$ for all coatoms $i$ in it, then it is a block 
ending in a 1; if not, then it is part of the final block, and the unique 
coatom with 
$|X_i| > 1$ is the uppermost element of the set, with GJW number equal to its size. 
We first consider the final block, and show how to deduce the corresponding part of 
the GJW sequence from it.

The final block consists of a number of sub-blocks, each of which is a connected 
component with a unique push move coatom, which is the highest row of the 
sub-block. We need to first determine the order of the sub-blocks, and then 
determine the GJW numbers and order of each sub-block.

\begin{prop}\label{higherpush}
Let $i$ and $x$ be push move coatoms, and let $r$ be the smallest element of $X_i$. Then $r$ covers 1 element in $X_{i, x}$ if $i<x$, and two elements otherwise.
\end{prop}
\begin{proof}
Suppose that $i<x$. There are no rooks above and to the left of either 
rook $(i, j)$ or rook $(x, y)$ in the initial placement; in particular, 
$y>j$. Rooks below row $i$ cannot be moved if one is to remain in $X_{i, 
x}$; similarly, rooks between rows $i$ and $x$ to the right of column $j$ 
cannot be moved, and rooks above row $x$ to the right of column $y$ cannot 
be moved. So the only rooks which can be moved in the placements in $X_{i, 
x}$ are the $i$-th and $x$-th rooks. Therefore, if $r-1$ is the number of 
empty columns to the left of $j$, and $s-2$ is the number of empty columns 
to the left of $y$, then $X_{i, x}$ is isomorphic to the rook poset of the 
partition $(r, s)$, with $2\le r<s$. The minimal element of $X_i$ 
corresponds to the rook placement $[1, s]$, which covers only $[1, s-1]$, 
while the minimal element of $X_x$ corresponds to the rook placement $[r, 
1]$, which covers both $[1, r]$ and $[r-1, 1]$. This completes the proof.
\end{proof}

This gives us the order of the push move coatoms, and consequently the order of the 
sub-blocks, each of which is a contiguous set of rows with highest row equal to its 
push move. Each sub-block corresponds to a tree in the ordered entanglement graph, 
whose root is the push coatom, identifiable as the only one with nontrivial 
$|X_i|$.

What remains is to deduce the ordering of the entanglement tree of each sub-block 
from the rook poset. The coatoms connected to the root each have a subtree below 
them, and each subtree is a set of contiguous rows. We need to figure out how to 
order these subtrees, which is tantamount to figuring out how to order the coatoms. 
Suppose we have two coatoms $j$ and $k$ connected to the push coatom $i$, and we 
are trying to figure out which one is higher.  As in the proof of 
Proposition~\ref{higherpush}, $X_{\{i, j, k\}}$ is isomorphic to the poset of those 
three rows, since moving any other rook creates a position below some other 
element. This poset is the rook poset of the partition $(s, s+1, s+1)$, where $s\ge 
4$ because the GJW number of the highest coatom is at least 2.

But now we can distinguish the higher of coatoms $j$ and $k$. Suppose 
$j<k$, so that $j$ corresponds to the first coatom in the rook poset 
of $(s, s+1, s+1)$, and 
$k$ to the second. The maximal rook placement is $[s, s+1, s-1]$. The rook 
placement for coatom $j$ is $[s-1, s+1, s]$, which is the maximal rook 
placement of $(s-1, s+1, s+1)$. The rook placement for coatom $k$ is $[s, 
s-1, s+1]$, which is the maximal rook placement of $(s, s, s+1)$. Since 
$s\ge 4$, these diagrams are not conjugate, and thus by induction on the 
number of boxes in the diagram we can tell apart their rook posets. So the 
lower coatom is the one whose lower interval in the rook poset is 
isomorphic to that of the partition $(s-1, s+1, s+1)$.



The same procedure works to tell apart children of elements which are 
not the push coatom. Suppose we are trying to figure out which of $j<k$ is higher, 
with both $j$ and $k$ being children of $i_n$, which is connected to the root $i_1$ 
via $i_{n-1}, \ldots, i_2$. Isolating the relevant rows as before, the rook poset 
$X_{\{j, k, i_n, \ldots, i_1\}}$ is isomorphic to the rook poset of $(s, s+1, s+1, 
\ldots, s+1)$, where $s \ge n+3$. 

The maximal element of this rook poset is $(s, s+1, s-1, s-2, s-3, \ldots, 
1)$. Then the rook placement $j$ is $[s-1, s+1, s, s-2, s-3, \ldots, 1]$, 
which is the maximal element of the rook poset of $(s-1, s+1, s+1, s+1 
\ldots, s+1)$, while the rook placement $k$ is $[s, s-1, s+1, s-2, s-3, 
\ldots, 1]$, which is the maximal element of the rook poset of $(s, s, 
s+1, s+1, \ldots, s+1)$. Again by the restrictions on $s$, these diagrams 
are not conjugate, so we can distinguish which is lower among the two 
coatoms.



\begin{figure}
\begin{center}
\includegraphics{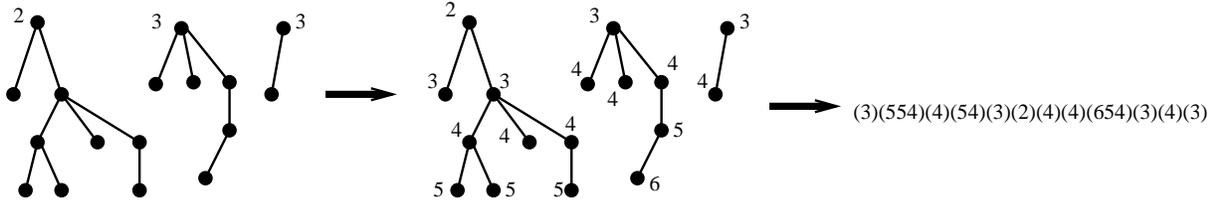}
\end{center}
\caption{\label{gjwrecover}
Starting with the ordered entanglement graph of the last block, along with 
the size of the $X_i$'s corresponding to the push coatoms, we can recover 
the GJW sequence. Using 
Proposition~\ref{whenswitch}, we can recover the GJW numbers of each 
coatom in the entanglement graph, and the order is given by the ordering 
of the entanglement graph. We have inserted parentheses grouping subtrees 
in the final GJW sequence for clarity.}
\end{figure}

This information allows us to recover the ordering of this portion of the entanglement graph. Once 
we have this, it is easy to recover the GJW sequence of the final block. We know the order and size 
of the sub-blocks, and for each sub-block the GJW number of the highest element $i$ is simply the 
size of $X_i$. We know exactly how the vertices of the sub-block tree (in the ordered entanglement 
graph) correspond to the other coatoms in the sub-block by the previous reasoning; given 
Proposition~\ref{whenswitch}, this allows us to reconstruct the entire GJW 
sequence; see Figure~\ref{gjwrecover}.

In this fashion, we have reconstructed the last part of the GJW sequence 
(the part which does not end in a 1) from the rook poset. We have 
identified the other blocks, and Theorem~\ref{mainthm} says that we cannot 
tell apart the order of the other blocks, so all that remains is to 
determine the structure of a given block ending in 1, up to conjugation, 
from its rook poset (which is $X_I$ in the original rook poset, $I$ being 
the collection of coatoms in the block.) We start by describing what this 
conjugation corresponds to in terms of the ordered entanglement graph. 

\begin{prop}
Given a block ending in 1, the ordered entanglement graph of its conjugate is the vertical 
reflection of its ordered entanglement graph (which is a tree).
\end{prop}

\begin{proof}

Augment the block by adding a 0 row. Then the ordered entanglement tree of a block is drawn by 
taking the maximal rook placement $\hat{1}$ and connecting each rook to the unique rook above it 
and to its left (the rook which it is switched with), except for the top coatom's connection, and 
turning this 45 degrees clockwise. The maximal rook placement $\hat{1}'$ of the conjugate block is 
just the conjugate of $\hat{1}$; this corresponds to reflecting the diagram, including rooks, about 
the upper-left to lower-right diagonal. After the rotation of 45 degrees is performed, this 
corresponds to reflecting the tree along a vertical line.

\end{proof}

Once we have the ordered entanglement tree, we can recover the GJW 
sequence of the block from the criterion given by 
Proposition~\ref{whenswitch}. We already know how to obtain the 
entanglement tree from the poset. All that remains is to figure out where 
the root of the tree (the top row) is, and figure out up to vertical 
reflection the order of the children of each node.

\begin{prop}\label{rootdist}

For each coatom $i$ in a given block, let $\hat{i}$ be an element of the 
rook poset with minimal rank which does not lie beneath $i$. Then if $i$ 
is a leaf, $r(\hat{i})$ is equal to the distance between $i$ and the root 
of the entanglement tree of that block, plus one (i.e. if a coatom is a 
child of the root, then it is 2.)

\end{prop}

\begin{proof}

Our goal is to find an element of minimal rank which is not less than $i$. 
By Lemma~\ref{belowcoat}, this is tantamount to saying that we have 
$\hat{i}$ of minimal rank such that $\{\hat{i}_1, \ldots, 
\hat{i}_i\}\not\leq \{i_1, \ldots, i_i\}$. Since $i$ is a leaf, there 
exists no rook below and to the right of rook $i$ in $\hat{1}$, so none of 
the first $i-1$ entries of $\lambda$ can be $\hat{1}_i$ or bigger, which 
implies that $\hat{1}_i = \lambda_i$. We claim that the rook 
placement $[1, 2, \ldots, i-1, 
\lambda_i, i, i+1, \ldots, \lambda_i-1, \lambda_i+1, i_i+2, \ldots]$ is 
such an element $\hat{i}$, i.e. of minimal rank not lying below $i$. 
Clearly it does not lie below $i$, since none of $i_1, \ldots, i_{i-1}$ 
can be bigger than $\lambda_i$, and $i_i$ is less than $\hat{1}_i = 
\lambda_i$. Furthermore, we claim that we picked this element by taking 
the minimal possible rook in each row to obtain an element less than $i$. 
This is clear for all rows except the $i$-th. For the $i$-th row, it 
follows easily from the fact that move $i$ moves rook $i$ to the rightmost 
possible position except for $\lambda_i$ and the positions occupied by the 
first $i-1$ rooks.

How many inversions are there in this rook placement? The only inversions come from $\lambda_i$ being bigger than some later element. This happens a total of $\lambda_i-i$ times, i.e. the $i$-th element in the GJW sequence of $\lambda$. However, the parent of each coatom is the rook with which it is switched, which is the first rook with smaller GJW number by Proposition~\ref{whenswitch}. This GJW number must be one less than the GJW number of the child, since GJW numbers can only decrease by 1 at any step. Proceeding in this fashion, the number of steps in the entanglement graph before we reach the root, which has GJW number 1, from coatom $i$, is precisely one less than the GJW number of $i$, which is $\lambda_i-i-1$. This completes the proof of the proposition.
\end{proof}

An easy induction shows that knowing how far each leaf is from the root of a tree allows us to 
determine the identity of this root. Drawing the tree with this root at the top, all that remains 
is to order the branches up to mirror reversal. This is tantamount to doing two things: first, at 
the highest branch point in the tree, for any three branches we need to be able to determine which 
is the middle one. Second, for any lower branch point, we need to be able to determine exactly 
which order the branches are in. We tackle these tasks in order.

Suppose we have three branches at the highest branch point in the tree. We 
can isolate the rows corresponding to those three children and to the 
chain of command heading back to the root, reducing us to the case where 
the diagram is $(s-2, s-1, s, s, \ldots, s)$, where there are $s$ rows in 
total and $s\ge 5$; the maximal element of the rook poset is $[s-2, s-1, 
s, s-3, s-4, \ldots, 1]$. We need to tell apart coatom 2, which is the 
rook placement $[s-2, s-3, s, s-1, s-4, \ldots, 1]$, from coatom 1, which 
is $[s-3, s-1, s, s-2, s-4, \ldots, 1]$, and coatom 3, which is $[s-2, 
s-1, s-3, s, s-4, \ldots, 1]$.  However, coatom 2 is the maximal element 
of the rook poset of the partition $(s-2, s-2, s, s, \ldots, s)$, while 
coatom 1 is the maximal element of the rook poset of $(s-3, s-1, s, s, 
\ldots, s)$, and 
coatom 3 is the maximal element of the rook poset $(s-2, s-1, s-1, s, s, 
\ldots, s)$. The 
last two boards are conjugate, but the board associated to coatom 2 is not 
conjugate to either of them, so again by induction on the number of the 
boxes in the diagram we can distinguish it from the others.

So, at the first branch point, given any three branches, we can determine which branch is the middle branch. 
What remains is, for any branch point which is not the first, to determine what order the branches are in, using 
information about the order of the branches at the higher levels of the tree. In other words, given two children 
of a common parent and a branching higher up the tree, we need to determine which of the children is higher. 
Without loss of generality (just conjugate the diagram), we can assume 
that the branch higher up the tree is to 
the right of the main branch we are considering. 

Again isolating the poset $X_R$ for the relevant set of rows, we are left 
with a copy of the rook poset of the Ferrers board $(s, s+1, s+1, \ldots, 
s+1, s+2, s+2, \ldots, s+2)$, where there are $k+1$ copies of $s+1$ and 
$s-k$ copies of $s+2$, for a total of $s+2$ entries. By the conditions 
given, we have $k\ge 1$, $s\ge 6$, and $s-k\ge 3$. The element $\hat{1}$ 
is $[s, s+1, s-1, s-2, \ldots, s-k, s+2, s-k-1, s-k-2, \ldots, 1]$. Thus 
rook $k+3$ is the branch higher in the tree. Our goal is to tell apart 
coatoms 1 and 2.

We use the previous technique of induction on the number of boxes in the 
diagram once again. Coatom 1 corresponds to the rook placement which 
replaces the first three rooks in $\hat{1}$ by $[s-1, s+1, s]$; this 
is the maximal element of the rook poset of the partition $(s-1, s+1, s+1, 
\ldots, s+1, s+2, \ldots, s+2)$, where there are still $k+1$ copies of 
$s+1$ and $s-k$ copies of $s+2$. Similarly, coatom 2 corresponds to the 
rook placement replacing these rooks by $[s, s-1, s+1]$, so its lower 
interval is isomorphic to the rook poset of the partition $(s, s, s+1, 
\ldots, s+1, s+2, \ldots, s+2)$, with $k$ copies of $s+1$ and $s-k$ copies 
of $s+2$. Again by the conditions on $s, k$, and $s-k$, these diagrams are 
not conjugate, so we can tell their rook posets apart, and consequently we 
can distinguish coatom 1 from coatom 2.


This completes the proof of Theorem~\ref{mainthm}. We recapitulate the algorithm for determining the GJW sequence 
from the rook poset.

\begin{enumerate}
\item Draw the entanglement graph of the coatoms. This is a forest by Corollary~\ref{tree}.
\item For each coatom $i$, determine the size of $X_i$.
\item All connected components with any coatom with $|X_i|>2$ are part of the final block by Proposition~\ref{whenpush}.
\begin{enumerate}
\item Each connected component consists of a number of contiguous rows in the final block. Proposition~\ref{higherpush} tells us the order of these sub-blocks.
\item For each sub-block, we can recover the root of the tree as the 
unique coatom with $|X_i|>2$, and we can recover the order of the branches 
at each level.
\item The last element of the GJW sequence of a sub-block is the size of $X_i$ for the top coatom, and we can 
recover the other elements in that sub-block by knowing the ordered tree and by Proposition~\ref{whenswitch}.
\end{enumerate}
\item Any other connected component is a block of its own.
\begin{enumerate}
\item We can determine the top coatom by Proposition~\ref{rootdist}. 
\item At the first branch point, we can determine the order of the branches up to mirror-reversal.
\item At any other branch point, we can determine the order of the branches by using the order of the branches further up the tree. Therefore, we can determine the entire rooted tree up to mirror-reversal.
\item The last element of the GJW sequence of a block is 1, and we can recover the other elements in that block by 
knowing the ordered tree
and by Proposition~\ref{whenswitch}. Mirror-reversing the tree corresponds to taking the conjugate Ferrers board.
\end{enumerate}
\end{enumerate}

\section{Applications and further questions} The two main applications of 
the resolution of this combinatorial problem were mentioned in the 
introduction. First, the category of rook posets of diagrams includes the 
category of lower intervals of 312-avoiding permutations in the Bruhat 
order on the symmetric group $S_n$. Therefore, Theorem~\ref{mainthm} has 
as a corollary a 
classification of poset isomorphism types of lower intervals of such 
permutations in the Bruhat order. Breaking up the set into blocks is the 
same as giving a nonzero composition of $[n]$ and considering the 
permutation which is the natural product of permutations on the components 
of the composition. Note that conjugating the diagram corresponds to 
taking the inverse permutation.

The steps used to prove Theorem~\ref{mainthm} have implications for the 
structure of these lower intervals. Proposition~\ref{whenpush} and 
Proposition~\ref{rootdist} are interesting statements about the Bruhat 
order, and the investigation of the $X_I$'s is an investigation of how the 
various subintervals fit together. In addition to classifying when 
312-avoiding permutations in the Bruhat order are isomorphic, 
Theorem~\ref{mainthm} also covers some non-312-avoiding permutations 
(permutations produced when the last part is greater than 1 may contain 
312-subsequences.) An obvious question to ask is the following 
generalization of the Bruhat translation of Theorem~\ref{mainthm}:

\begin{question}
When are two lower intervals in the Bruhat order isomorphic?
\end{question}

The other main application is to Ding's Schubert varieties. Schubert cells are in bijection with rook placements, 
with the closure of a Schubert cell equal to the union of the cells of all rook placements which lie below it. 
Proposition~\ref{whenpush} tells us about the occurrence of simple substructures (self-contained chains of cells), 
and Theorem~\ref{mainthm} gives us an exact categorization of when the poset cell structures of two Ding varieties 
are the same.

This equivalence fits in between two previous notions of equivalence. 
Ordinary rook equivalence corresponds to two Schubert varieties having 
identical Poincar\'{e} series~\cite{Ding2}, or to their integral 
cohomology rings 
being isomorphic as graded Abelian groups. The finer rook equivalence 
given 
by Theorem~\ref{DMR-thm} tells us when two Ding varieties are isomorphic, 
or when their cohomology rings are isomorphic as graded rings. The 
condition of Theorem~\ref{mainthm} is, from the combinatorics, visibly 
stronger than the first and weaker than the second.

Comparing Theorem~\ref{mainthm} to Theorem~\ref{DMR-thm} tells us when two 
Schubert varieties can have the same poset of Schubert cells but not be 
isomorphic or have identical cohomology. These situations are generated by 
three operations: adding 0's after 1's, removing 0's, and conjugating 
blocks ending in 1's (one can think of this as adding a 0, conjugating a 
block ending in 0, and then removing the 0.) 

\begin{question}
What does the equivalence condition of Theorem~\ref{mainthm} 
tell us about the cohomology rings of the corresponding Schubert 
varieties? In other words, do these operations preserve 
some interesting invariant of the cohomology ring? What do the operations 
of conjugation and adding 0's after 1's do to the Schubert variety?
\end{question}

Adding 0's after 1's changes the Schubert variety from a nontrivial fiber 
bundle to a straightforward product; can we say more? What does 
conjugation do to the Schubert variety?

\section{Acknowledgements}
I would like to thank Vic Reiner for comments on an early draft of this 
paper, and Jeremy Martin for observing the first sentence of the last 
paragraph.


\begin{thebibliography}{99}

\bibitem{DMR}
M.~Develin, J.L.~Martin, and V.~Reiner, Classification of Ding's Schubert 
varieties: finer rook equivalence, 
to 
appear in {\it Canadian J. Math}, \texttt{http://arxiv.org/math.AG/0403540}.

\bibitem{Ding}
K.\ Ding,
Rook placements and cellular decompositions of partition varieties.
{\it Discrete Math.} {\bf 170} (1997), 107-151.

\bibitem{Ding2}
K.\ Ding,
Rook placements and classification of partition varieties
$B\backslash M\sb \lambda$,
{\it Commun.\ Contemp.\ Math.} {\bf 3} (2001), 495--500.

\bibitem{FS}
D.\ Foata and M.-P.\ Sch\"utzenberger,
On the rook polynomials of Ferrers relations.
Combinatorial theory and its applications, II (Proc.\ Colloq., Balatonf\"ured, 1969),
pp.~413--436. North-Holland, Amsterdam, 1970.

\bibitem{GR}
A.M.\ Garsia and J.B.\ Remmel,
$Q$-counting rook configurations and a formula of Frobenius.
{\it J.\ Combin.\ Theory Ser.\ A} {\bf 41} (1986), 246--275.

\bibitem{GJW}
J.R.\ Goldman, J.T.\ Joichi, and D.E.\ White,
Rook theory. I. Rook equivalence of Ferrers boards.
{\it Proc.\ Amer.\ Math.\ Soc.} {\bf 52} (1975), 485--492.

\bibitem{Kapl}
I.\ Kaplansky and J.\ Riordan,
The problem of the rooks and its applications.
{\it Duke Math.\ J.} {\bf 13} (1946), 259--268.

\bibitem{MRW}
J. Martin, V. Reiner, and J. Wagner, personal communication.

\end{thebibliography}
\end{document}